\documentclass[smallextended]{svjour3}       
\smartqed  
\usepackage{graphicx}
\usepackage{subfigure}
\usepackage{mathptmx}      
\usepackage{amsfonts}
\usepackage{color}
\usepackage[ruled]{algorithm2e}
\usepackage{epstopdf}
\usepackage{epsfig}

\begin{document}

\title{Adaptive Block Coordinate DIRECT algorithm
}


\author{Qinghua Tao        \and
        Xiaolin Huang     \and
        Shuning Wang  \and
        Li Li(Corresponding Author)
}


\institute{Qinghua Tao \at
              Department of Automation, Tsinghua University, Beijing 100084, PR China \\
              \email{taoqh14@mails.tsinghua.edu.cn}           
          \and
          Xiaolin Huang  \at
          Institute of Image Processing and Pattern Recognition, Shanghai Jiao Tong University, Shanghai 200240, PR China\\
        \email{xiaolinhuang@sjtu.edu.cn}
           \and
        Shuning Wang \at
          Department of Automation, Tsinghua University, Beijing 100084, PR China \\
        \email{swang@mail.tsinghua.edu.cn}
       \and
           Li Li(Corresponding Author)  \at
             Department of Automation, Tsinghua University, Beijing 100084, PR China \\
            \email{li-li@mail.tsinghua.edu.cn}\\
            \\
            This project is jointly supported by the National Natural Science Foundation of
China (61473165, 61134012, 61603248) and the National Basic Research Program of China (2012CB720505).
}

\date{Received: date / Accepted: date}

\maketitle

\begin{abstract}

 DIviding RECTangles (DIRECT) is an efficient and popular method in dealing with bound constrained optimization problems. However, DIRECT suffers from dimension curse, since its computational complexity soars when dimension increases. Besides, DIRECT also converges slowly when the objective function is flat. In this paper, we propose a coordinate DIRECT algorithm, which coincides with the spirits of other coordinate update algorithms. We transform the original problem into a series of sub-problems, where only one or several coordinates are  selected to optimize and the rest keeps fixed. For each sub-problem, coordinately dividing the feasible domain enjoys low computational burden. Besides, we develop adaptive schemes to keep the efficiency and flexibility to tackle different functions.  Specifically, we use block coordinate update, of which the size could be adaptively selected, and we also employ sequential quadratic programming (SQP) to conduct the local search to efficiently accelerate the convergence even when the objective function is flat. With these techniques, the proposed algorithm achieves promising performance on both efficiency and accuracy in numerical experiments.

\end{abstract}

\section{Introduction}
\label{intro}
Bound constrains optimization problems have been applied in many applications, such as operation research, engineering, biology and biomathematics, chemistry, finance, etc \cite{Ref1,Ref11,Ref12,Ref13,Ref14,Ref15}.  Mathematically, this kind of problems can be formulated as the following,
\begin{equation}\label{general_problem}
\min_{x\in \Omega} f(\mathbf{x}),
\end{equation}
where $\Omega = \{\mathbf{x }\in \mathbb{R}^n: l_i \leq x_i \leq u_i, i=1,...,n\}$, $\mathbf{l},\mathbf{u}\in \mathbb{R}^n$ and $f: \mathbb{R}^n \rightarrow \mathbb{R}$.
Notice that the convexity on $f(\mathbf{x})$ is not required, and thus global optimization techniques are needed.

There have been many efficient and important methods for global optimization, such as branch-and-bound methods and some stochastic methods \cite{Ref1,Ref4,Ref5}. In recent years, DIviding RECTangles (DIRECT), based on the strategy of branch-and-cut method, has shown very good performance in dealing with problem (\ref{general_problem}), i.e., problems with bound constraints, since the basic idea of DIRECT is to do domain partitions \cite{Ref11,Ref13,Ref8,Ref9,Ref10}.  In low dimensions, domain partition is quite effective and hence DIRECT has obtained success in some applications \cite{Ref1,Ref11,Ref13,Ref14,Ref15}.

The existing DIRECT algorithms consider the whole domain and its computational complexity soars exponentially with dimension increasing. DIRECT looses its efficiency on both computational time and required storage space for high-dimensional problems. In fact, this phenomenon dose not only appear in domain-partition-based methods, but also happens in gradient-based methods. For the latter, the combination with coordinate descent method (CDM) has become a promising alternative and has been widely applied in image reconstruction \cite{Ref37}, dynamic programming \cite{Ref38}, flow routing \cite{Ref39} and the dual of a linearly
constrained \cite{Ref40,Ref41}, etc. A similar idea that coordinately partitions the domain is also applicable for domain-partition-based methods. Therefore, in this paper, we propose a coordinate DIRECT algorithm. To be specific, we transform the original problem into many sub-problems, where we only select one coordinate to conduct DIRECT and keep the rest fixed. Low computational burden can be expected with coordinate update, but we need to guarantee that there is efficient descent in each sub-problem for different objective functions. Hence, we introduce a switch to the combination with block coordinate descent method (BCDM) when necessary. We also employ SQP to locally find a good solution. Since SQP is efficient when the objective function performs smooth, the employment of SQP can accelerate the descent when the objective function is flat, where original DIRECT methods are reported to loose efficiency in convergence and descent \cite{Ref1,Ref4,Ref9,Ref10}. Summarizing the above discussion, we establish the adaptive block coordinate DIRECT (ABCD) algorithm: its basic idea is coordinate update and it can adaptively choose single coordinate partition, block-coordinate partition, and local optimizer SQP. In numerical experiments, the proposed ABCD algorithm presents very promising performance in comparisons with DIRECT and other relevant algorithms.

The rest of this paper is organized as follows. Section 2 gives a review on the background and preliminaries. Section 3 establishes the proposed algorithm, i.e., ABCD algorithm. A series of numerical experiments are conducted to demonstrate the overall performance of the proposed ABCD algorithm in Section 4. Section 5 ends this paper with conclusions.

\section{Backgrounds}
DIRECT was first introduced by Jones \cite{Ref8,Ref9} to search the global optimum of a real valued objective function with bounded constraints, i.e., problem (1). DIRECT emerges as a natural extension and generalization of the Lipschitz optimization (LO) methods proposed by Pijavskiy and Shubert \cite{Ref17}. For problem (1), it is customary to make assumptions on the objective function $f(\mathbf{x})$ for continuity, Lipschitz continuity, smoothness or differentiability \cite{Ref3}. LO requires Lipschitz continuity in searching domain, and it also demands the knowledge of the Lipschitz constant. These requirements for the objective function are sometimes too demanding, which narrows the applications of LO. Different from LO, DIRECT needs no knowledge of Lipschitz constant, and DIRECT even does not require the objective function to be Lipschitz continuous\cite{Ref8}. It only requires the objective function to be continuous in the neighborhood of the global optimum, which enables it to deal with a wider range of problems.

Considering problem (1), DIRECT starts from unifying the searching domain $\Omega$ into a unit hyper-cube $\bar{\Omega}$. That is
\begin{equation}
\bar{\Omega}={\mathbf{x}\in \mathbb{R}^n:0\leq x_i \leq 1}, i=1,...,n.
\end{equation}
 DIRECT optimizes the problem over the normalized hyper-cube, of which the center is denoted as $\mathbf{c}_1$. In the initialization, DIRECT samples the center $\mathbf{c}_1$, then we have the objective function $f(\mathbf{c}_1)$, which is denoted as the initial optimum $\hat{f}^*$.

Next, DIRECT turns to divide the hyper-cube. DIRECT evaluates the objective function $f(\mathbf{x})$ at points $\mathbf{c}_1 \pm \delta \mathbf{e}_i, i=1,...,n$, where $\delta$ is one third of the side length of the hyper-cube and $\mathbf{e}_i$ is the unit vector of $i$th dimension, i.e., a vector with a one in the $i$th position and zeros elsewhere. An easy way to divide the hyper-cube is selecting one dimension arbitrarily and splitting it along this dimension. However, arbitrariness is not efficient and DIRECT heuristically uses the following criterion:
\begin{equation}
w_i = \min (f(\mathbf{c}_1+\delta \mathbb{e}_i),f(\mathbf{c}_1 - \delta \mathbb{e}_i)), i=1,...,n,
\end{equation}
and divides the hyper-cube starting with the smallest $w_i$ into thirds. Then $\mathbf{c}_1 \pm \delta \mathbb{e}_i$ is the center of new hyper-rectangle.

After the division, DIRECT proceeds to the identification of potentially optimal hyper-rectangles (POHs), which have potential to obtain better solutions, even the global optimum \cite{Ref8}. In DIRECT, a constant $\widetilde{K}$ is introduced to determine POHs. The basic idea of DIRECT is to explore better solutions among all the POHs. More precisely, DIRECT samples all ``potentially optimal" hyper-rectangles as defined below, and this pattern is repeated  in each iteration \cite{Ref8}.

\textbf{Definition 1} Let $\epsilon$ be a positive constant and $f_{\min}$ be the current minimum. Interval $j$ is said to be potentially optimal if there exists a rate-of-change constant $\widetilde{K}$ such that
\begin{equation}
\begin{array}{lll}
f(\mathbf{c}_j) - \widetilde{K} \mathbf{d}_j & \leq & f(\mathbf{c}_i) -  \widetilde{K} \mathbf{d}_j,~\forall i\\
f(\mathbf{c}_j) - \widetilde{K} \mathbf{d}_j  & \leq  & f_{\min} - \epsilon |f_{\min}|,
\end{array}
\end{equation}
where $\mathbf{d}_j$ is a measure for hyper-rectangle $j$. Researchers commonly choose the distance from center $\mathbf{c}_j$ to its vertices as the measure \cite{Ref8,Ref18}. The first condition in the definition forces the POHs to obtain the lowest objective values among the rectangles which share the same measure $\mathbf{d}_j$. The second condition insists that the POHs, based on the rate-of-change constant $\widetilde{K}$, exceed the current best solution by a nontrivial amount \cite{Ref8}. Parameter $\epsilon$ is usually chosen as $10^{-4}$ by researchers, since the experimental data show that $\epsilon$ has a negligible effect on the calculation when it is chosen within $10^{-2}$ to $10^{-7}$ \cite{Ref18}.

  Once a hyper-rectangle is identified to be the POH, then it needs to be divided into smaller ones. The dividing procedure is restrained to be conducted only along the dimensions with the longest side-length, which ensures a shrinkage in every dimension \cite{Ref18}. The sampling and dividing rules of DIRECT are shown in Algorithm 1. The formal description of DIRECT algorithm is presented in Algorithm 2.

\begin{algorithm}
\caption{Sampling and Dividing Algorithm \cite{Ref8} }
\footnotesize{
Initialization: Identify the set $I$ of dimensions with the maximum side length $d_s$, set $\delta = d_s/3$.\\
Sampling: Sample $f(x)$ at $\mathbf{c}\pm \delta \mathbf{e}_i, i\in I$, where $\mathbf{c}$ is the center point of the hyper-rectangle.\\
Dividing: Divide the hyper-rectangle into thirds alongside the dimensions in $I$, starting with the dimension with the smallest $w_i$, where $w_i=\min f(\mathbf{c}\pm \delta \mathbf{e}_i)$.
}
\end{algorithm}

\begin{algorithm}
\caption{DIRECT Algorithm \cite{Ref8} }
\footnotesize{
\KwIn{$f, \epsilon, N_1, N_2$}
\KwOut{$f_{\min}, \mathbf{x}_{\min}$}
Normalize the search space to be the unit hypercube with center point $\mathbf{c}_1$.\\
Evaluate $f(\mathbf{c}_i)$, $f_{\min} = f(\mathbf{c}_1)$.\\
Set the number of iterations $t=0$, and function evaluations $m=1$.\\
\While{$t<N_1$ and $m<N_2$}
{
Identifying the set $S$ of potentially optimal hyper-rectangles.\\
\While{$S\neq\emptyset$}
{
Take $j\in S$.\\
Sample new points, evaluate $f$ at the new points and divide the hyper-rectangles with Algorithm 1.\\
Update $f_{\min}$, $x_{\min}$ and $m=m+\triangle m$, where $\triangle m $ is the number of new points sampled.\\
Set $S=S - \{j\}$.
}
$t = t+1$.
} }
\end{algorithm}

In DIRECT, the progress of the optimization is governed only by the evaluations of the objective function. It iteratively divides the longest side of the selected POHs and obtains several sub-rectangles. There are two key processes required in every iteration of optimization. The first process is to identify the POHs, which are identified to  potentially contain good, unsampled points. DIRECT sorts the sub-rectangles by their measures $\mathbf{d}_j$, and selects the rectangles with the lowest function values at centers from every sorted measure $\mathbf{d}_j$. POHs are identified from  the set of these rectangles with Definition 1.  The second process of DIRECT is to sample and divide POHs, which shrinks the location of the global optimum into smaller space. The repetition of the two key processes can iteratively approach the global optimum and confine the global optimum within a very small rectangle. Therefore, DIRECT iteratively approaches  better solutions as iterations go on. Without the limit of iterations and function evaluations, DIRECT is guaranteed with global convergence and it can restrict the global optimum to a small rectangle with any given accuracy.

Unfortunately, the efficiency of DIRECT drops rapidly with dimension increasing, since working with all the coordinates of an optimization problem at each iteration may be inconvenient and the required function evaluations are soaring with dimensions increasing. From \cite{Ref8}, we know that, after $r$ divisions, the rectangle have $p=$ mod$(r,n)$ sides of length $3^{-(k+1)}$ and $n-p$ sides of length $3^{-k}$, where $k=(r-p)/n$. Thus, the center-to-vertex distance of the rectangle is given by
\begin{equation}\label{eq5}
d = 0.5\sqrt{[3^{-2(k+1)}p+3^{-2k}(n-p)]}.
\end{equation}

From equation (\ref{eq5}), it is easy to obtain the smallest required division $r$ of a rectangle with the center-to-vertex distance $d$. In DIRECT, a rectangle with the center-to-vertex distance $d$ is required to undergo at least $r = n\log_3(\sqrt{n}/2d)+p-n$ divisions. Since new rectangles are formed by dividing existing ones into thirds on the longest side, every division is accompanied with 2 function evaluations. Thus, in DIRECT, any rectangle with the center-to-vertex distance $d$ is accompanied with at least $N_0 = \lceil2^{n\log_3(\sqrt{n}/2d)+p-n}\rceil$ function evaluations. When $n=3$ and $\epsilon=10^{-4}$, $N_0 =\lceil 581.09\rceil$. When $n=10$ and $\epsilon=10^{-4}$, $N_0 =\lceil 7.3056^{10}\rceil$. We can see that the computational complexity of the DIRECT soars exponentially with dimension $n$ increasing.

 In high dimensions, DIRECT suffers from dimension curse, which is also a common defect of algorithms based on strategies of domain partition and function evaluation. Besides the effect of dimension, DIRECT also converges slowly around smooth area, where the objective function is very flat. In such case, even DIRECT gets close to the basin of global optimum, the smooth neighbor around the optimum may also hamper it to achieve higher accuracy.
 Along with the spirits of DIRECT,  multilevel coordinate search (MSC) is proposed  with the strategy of domain partition and function evaluation to tackle problem (1) \cite{Ref34}. MCS partitions the domain into small boxes with more irregular splitting. In contract to DIRECT, MCS simplifies the division procedure. To speed up the convergence, MCS introduces a local enhancement to start local searches when the corresponding sub-domains reach the maximal split level.  Although MCS improves computing speed compared with DIRECT, it lacks further exploration when all the local searches are finished. Thus the accuracy of MCS is undesirable if the local search fails to bring the solution to the sufficiently small neighborhood of the global optimum with given accuracy. Recently, simultaneous optimistic optimization (SOO) is introduced by Munos to solve problem (1) \cite{Ref3}. Similar to MCS, SOO is also closely related to DIRECT algorithm. It works by partitioning the domain into sub-parts, namely cells. SOO chooses to split the cell with the smallest value at its center.  As its name suggests, the basic idea of SOO is to optimally choose the sub-domain whose objective value is the lowest in the corresponding depth. The basic idea of SOO is very simple, but SOO is a global optimization algorithm which has a finite-time performance under some weak assumptions on the objective function\cite{Ref32}. However, SOO still suffers from dimension curse. Moreover, SOO is hard to approach the global optimum when the surface of the objective function is complicated or the problem has many local optima.  DIRECT, SOO and MCS all pose global convergence theoretically, and their global convergence comes when when it is allowed to have enough computational time and splitting depth \cite{Ref1,Ref8,Ref23}.

\section{Adaptive Block Coordinate DIRECT algorithm}
 DIRECT has shown to be efficient in low dimensions \cite{Ref1,Ref8,Ref23,Ref19,Ref20}, but its computation complexity soars exponentially with dimension increasing, and working with all the coordinates at each iteration is inconvenient.  Inspired by CDM, degenerating the original problems into many one-dimensional sub-problems is expected to bring fast speed. In this paper, we present the ABCD algorithm, whose basic idea is conducting the optimization coordinately to maintain the efficiency of DIRECT in low dimensions, so as to improve the speed. Instead of considering the optimization on all coordinates, we transform the original problem into a series of sub-problems, where we select only one coordinate to optimize with one-dimensional DIRECT. However, the coordinate update still possibly brings trivial descent in each sub-problem, then we modify it to be adaptive. In such case, we employ SQP to conduct a local search and then switch to the block coordinate update. The local optimizer helps to accelerate the convergence, while the switch explores further improvements with larger size of chosen coordinates, which brings more flexibilities.
\subsection{DIviding RECTangles on One Coordinate}
 DIRECT processes one-dimensional function very fast. In order to maintain DIRECT's efficiency in low dimensions and tackle high-dimensional problems as well, we introduce to develop a coordinate update DIRECT algorithm.

 The main concept is to optimize one single coordinate each time, and keep the rest fixed to constants. Coordinate update has been widely applied in the optimization with convex objective functions and proved to efficiently reduce computational complexity \cite{Refbca}.  This basic pattern is to decrease problem dimension from $n$ to 1, which formulates the original problem into many one-dimensional sub-problems. For high-dimensional problems, exponentially soaring function evaluations lead to slow convergence. Usually, DIRECT exhausts its max iterations or max dividing depth before reaching the given accuracy. By contrast, the degeneration to low dimensions sharply relives computational burden and enables CPU memory to get released at the end of every sub-problem. We conduct sampling and dividing on just one coordinate, then the solution can be quickly restricted to sufficiently small space, which enables solve high-dimensional problems with fast speed. Even in the cases where DIRECT can reach the given accuracy with desirable time consumption for some problems, the coordinate update helps DIRECT improve its speed to a higher level. Take Levy function for example. Levy function belongs to Hedar test set which is regarded as a benchmark for global optimization \cite{Ref1}. DIRECT solves Levy function from 6 dimension to 18 dimension with given accuracy $\epsilon=10^{-4}$, which is defined as the difference between the objective value obtained by the test algorithm and the truly global minimum. The results are shown in Fig.1. DIRECT takes only no more than 7 seconds to reach the given accuracy $\epsilon=10^{-4}$ for  Levy function from 6 dimensions to 18 dimensions.  With coordinate update, the running time for Levy function decreases to less than 0.2 second. Thus, the coordinate DIRECT not only remains to solve the problem with given accuracy, but also greatly increases the efficiency.

  \begin{figure}\label{levy_toy}
\centering
\includegraphics[width=0.6\columnwidth]{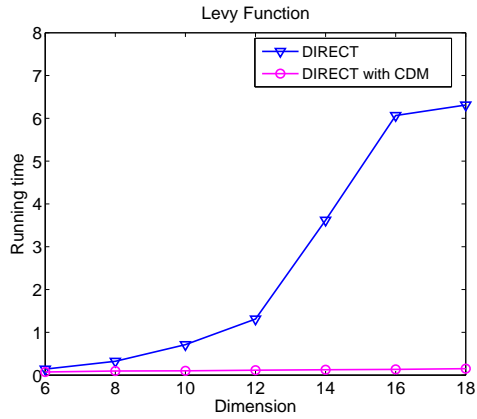}
\caption{Running time (s) of Levy function in dimension $n=6,8,...,18$ for reaching accuracy $\epsilon = 1e-4$. `DIRECT' means the results of directly applying DIRECT to Levy function. `DIRECT with CDM' represents the results of applying DIRECT with coordinate update to Levy function.  }
\end{figure}

Besides computational burden in high dimensions, DIRECT suffers slow convergence in a smooth area, which also affects the speed and accuracy. But when we update the problem coordinately, this problem can be overcome. To demonstrate this characteristic, we consider the function $f(x) = 0.1(x-0.4)^6$, which is very flat in interval $[-0.6, 1.4]$. We apply one-dimensional DIRECT to function $f(x) = 0.1(x-0.4)^6$ subjected to $x\in [-1, 1]$. The results are shown  in Fig.2.

 \begin{figure}\label{exp}
\centering
\subfigure[] { \label{fig:a}
\includegraphics[width=0.45\columnwidth]{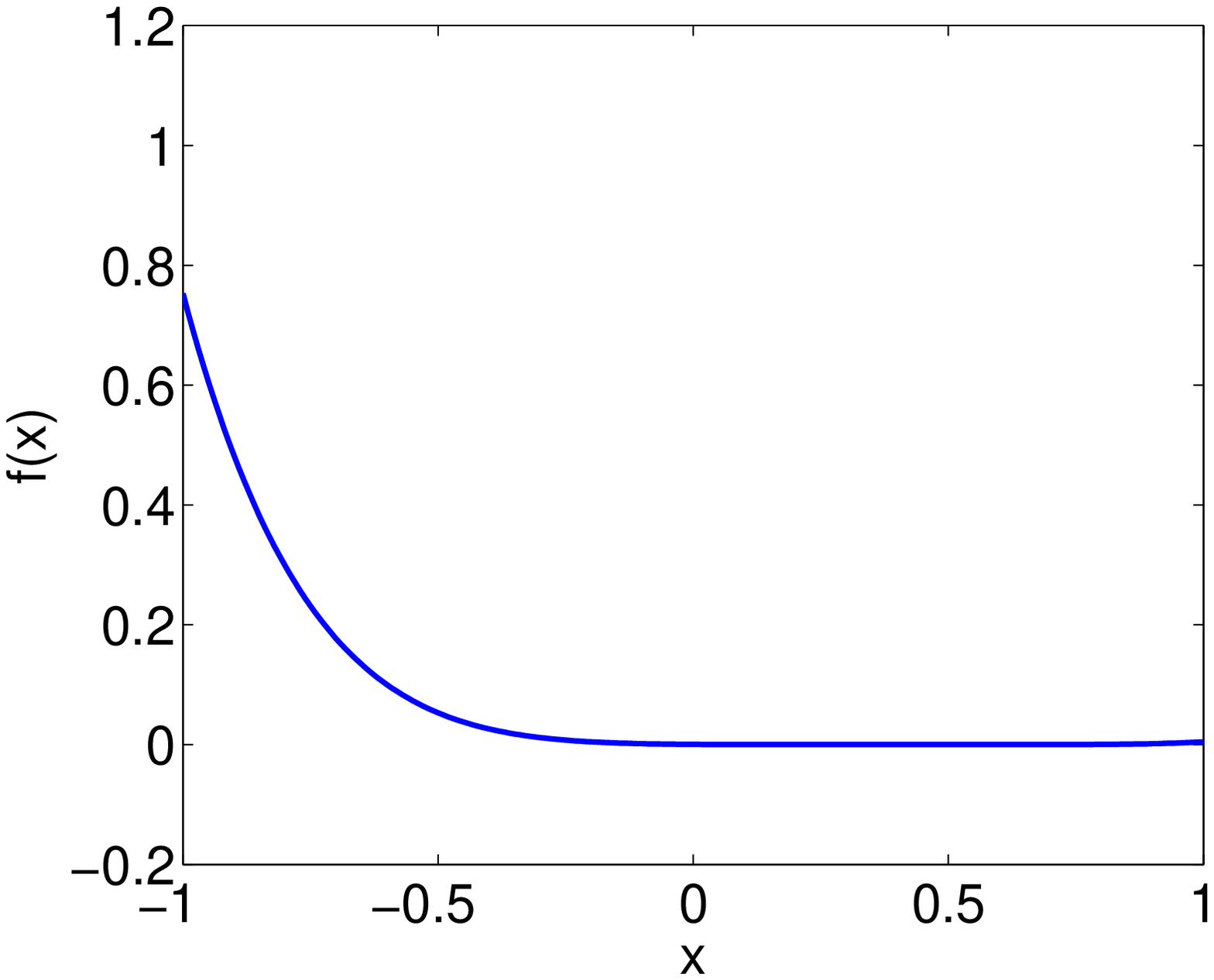}
}
\subfigure[] { \label{fig:b}
\includegraphics[width=0.45\columnwidth]{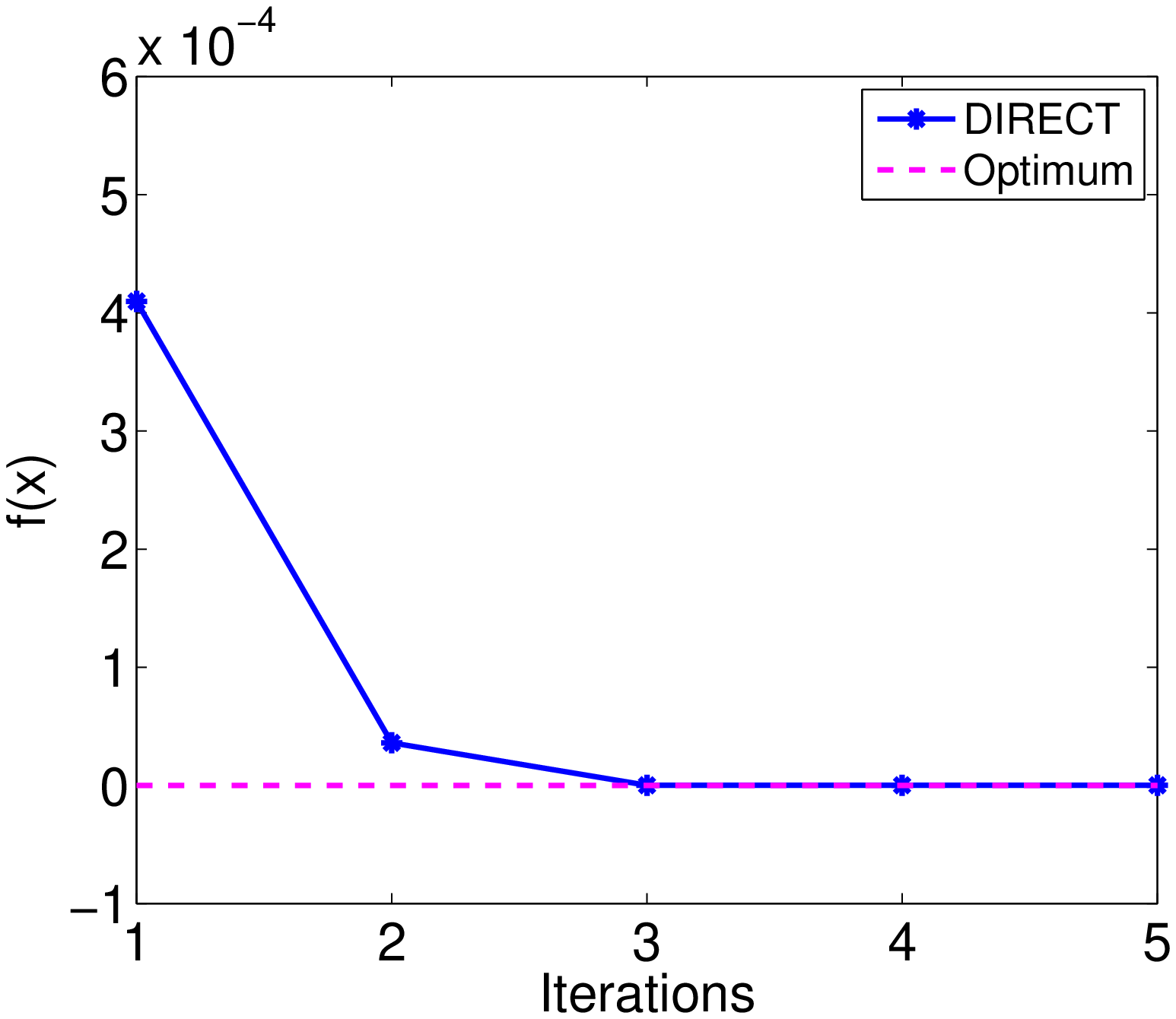}
}
\caption{Simulation results of function $f(x) = 0.1(x-0.4)^6$ when applying DIRECT. Subfigure (a) is the curve of function $f(x) = 0.1(x-0.4)^6$ in interval $[-1, 1]$. Subfigure (a) illustrates the results of DIRECT in every sub-problem. `Optimum' represents the global minima of function  $f(x) = 0.1(x-0.4)^6$ .   }
\end{figure}

It can be seen from Fig.2 that one-dimensional DIRECT works very well even when the objective function is very flat. For function $f(x) = 0.1(x-0.4)^6$, DIRECT only takes 5 iterations to achieve accuracy $\epsilon$. Even the degenerated sub-problem is very flat, the one-dimensional DIRECT can efficiently solve the sub-problem. Besides dimension degeneration, the fast convergence in flat objective function is another merit of coordinate update.

Instead of repeating sampling and dividing on the whole domain, we conduct DIRECT on just one coordinate. In the proposed algorithm, we first generate a starting point $\mathbf{x}_0$ in feasible domain $\Omega$, which satisfies the bounded constraints. In each sub-problem, we select the $i$th coordinate $\mathbf{x}(i)$ of $\mathbf{x}$ to optimize and keep the rest fixed to $\mathbf{x}_0(j),j=1,...,n,j\neq i$. Thus, the procedures of sampling, dividing, and identifying POHs are all conducted in the selected coordinate $\mathbf{x}(i)$. When the one-dimensional DIRECT algorithm converges on coordinate $\mathbf{x}(i)$, we sequentially select the adjacent coordinate $i=i+1$ to continue the next sub-problem. Algorithm stops when the objective function $f(\mathbf{x})$ no longer decreases among sub-problems or the optimum is restricted to a sufficiently small rectangle. This is the basic and simplest form of the proposed algorithm.

The original problem may possibly reach a saturation in descent with coordinate DIRECT, thus the degeneration to one dimension possibly lead to local optimum. Moreover, even algorithm runs very fast with the selected coordinate in each sub-problem, the descent of the objective in each sub-problem can still be quite small. In such case, if we continue selecting only one coordinate to optimize, it is possible that the original objective function $f(\mathbf{x})$ shares very sparse contours and each sub-problem only achieves a trivial descent in the objective. Although each sub-problem runs very fast in convergence, it requires tremendous numbers of sub-problems to reach the global optimum. To conquer these drawbacks, Section 3.2 and Section 3.3 are established.

\subsection{DIviding RECTangles on Block Coordinates }

The strategy of coordinate update can be extended. In every sub-problem, we can adaptivley select $1 < m\leq n$ coordinates to optimize and keep the rest fixed to constants, which resonates with the core idea of BCDM \cite{Ref29,Ref30}. When $m=n$, the proposed algorithm degenerates to DIRECT. Thus, the proposed algorithm can be regarded as an extension and generalization of DIRECT algorithm from this view.

 Differently from Section 3.1, we partition the coordinates into several blocks. In each sub-problem, we focus on updating a single block only, and keep the remaining blocks fixed. Sequentially optimizing each coordinate with DIRECT algorithm is the basic and the simplest form of the proposed algorithm, where each coordinate gets the same chance to be optimized. However, as mentioned in the end of Section 3.1, when the objective function $f(\mathbf{x})$ appears trivial descent in sub-problems and even reaches a local optimum, adaptively selecting more coordinates to do the optimization may possibly lead further descent and more efficiency to the objective function $f(\mathbf{x})$, since the objective function $f(\mathbf{x})$ may help bypass the local optimum and enjoy sharper descent in each sub-problem with larger size of chosen coordinates. Therefore, the proposed ABCD algorithm is established.

Take Rosenbrock function to illustrate this statement. Rosenbrock is regarded as a benchmark function in global optimization, and it is shown to be difficult to be tackled \cite{Ref1}. We try to solve the 12-dimensional Rosenbrock function with ABCD of coordinate update and blcok coordniate update. We terminate the algorithm in different running time, which ranges from 2s to 6s with the interval of 0.5s. Fig.3 shows the CPU running time against the absolute error.

\begin{figure}[!hbt]\label{rosen_bcdm_1}
\centering

\includegraphics[width=0.6\linewidth]{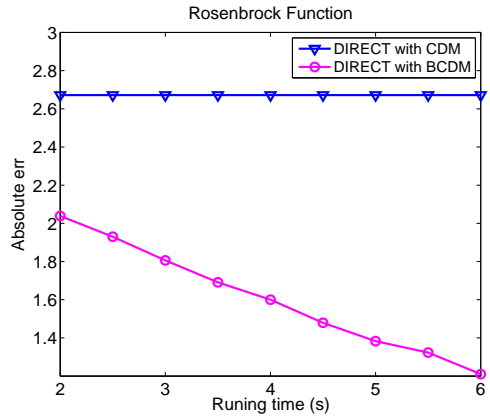}

\caption{Simulation results of optimizing Rosenbrock function. `DIRECT with CDM' represents the absolute error when sequentially selecting one coordinate to conduct one-dimensional in each sub-problem. `DIRECT with BCDM' denotes the absolute error when sequentially choosing two coordinate to optimize in each sub-problem.  }
\end{figure}

From Fig.3, we can see that sequentially conducting DIRECT on one coordinate leads to a local minimum within 2s. When algorithm runs longer, the coordinate update fails to get out of the current local minimum. However, sequentially conducting DIRECT on two coordinates bypasses the local minimum which the coordinate update gets stuck in. With more computational time, the block coordinate update gradually decreases the absolute error.  Therefore, when the coordinate update fails to achieve higher accuracy and converges inefficiently, we can switch to block forms. In the switch, we randomly select $m\geq 2$ coordinates to conduct DIRECT in each sub-problem. The switch is designed to check if there is any descent with other different size of chosen coordinates, which possibly helps check further improvements to approach a better solution, even the global optimum.

\subsection{Local optimizer }

 Although one-dimensional DIRECT works well in fast speed, it may appear trivial descent on the objective function $f(\mathbf{x})$ in each sub-problem. There are two possible conditions for this phenomenon. The first condition is that one-dimensional DIRECT fails to bring further improvements. The second one is that the original objective function $f(\mathbf{x})$ arrives in a smooth area, where $f(\mathbf{x})$ is very flat. For the first condition, as mentioned in Section 3.2, the switch to block coordinate update provides possibilities to conquer it. However, if it is the latter, the switch to block coordinate update may possibly remain failing in achieving sharp descent in each sub-problem. Take Rosenbrock function for example. In order to observe the contour map, we tackle the two-dimensional Rosenbrock function with coordinate update. Below, Fig.4 depicts the searching trace of the coordinate DIRECT.
\begin{figure}[!hbt]
\centering
\subfigure[] { \label{fig:a}
\includegraphics[width=0.45\linewidth]{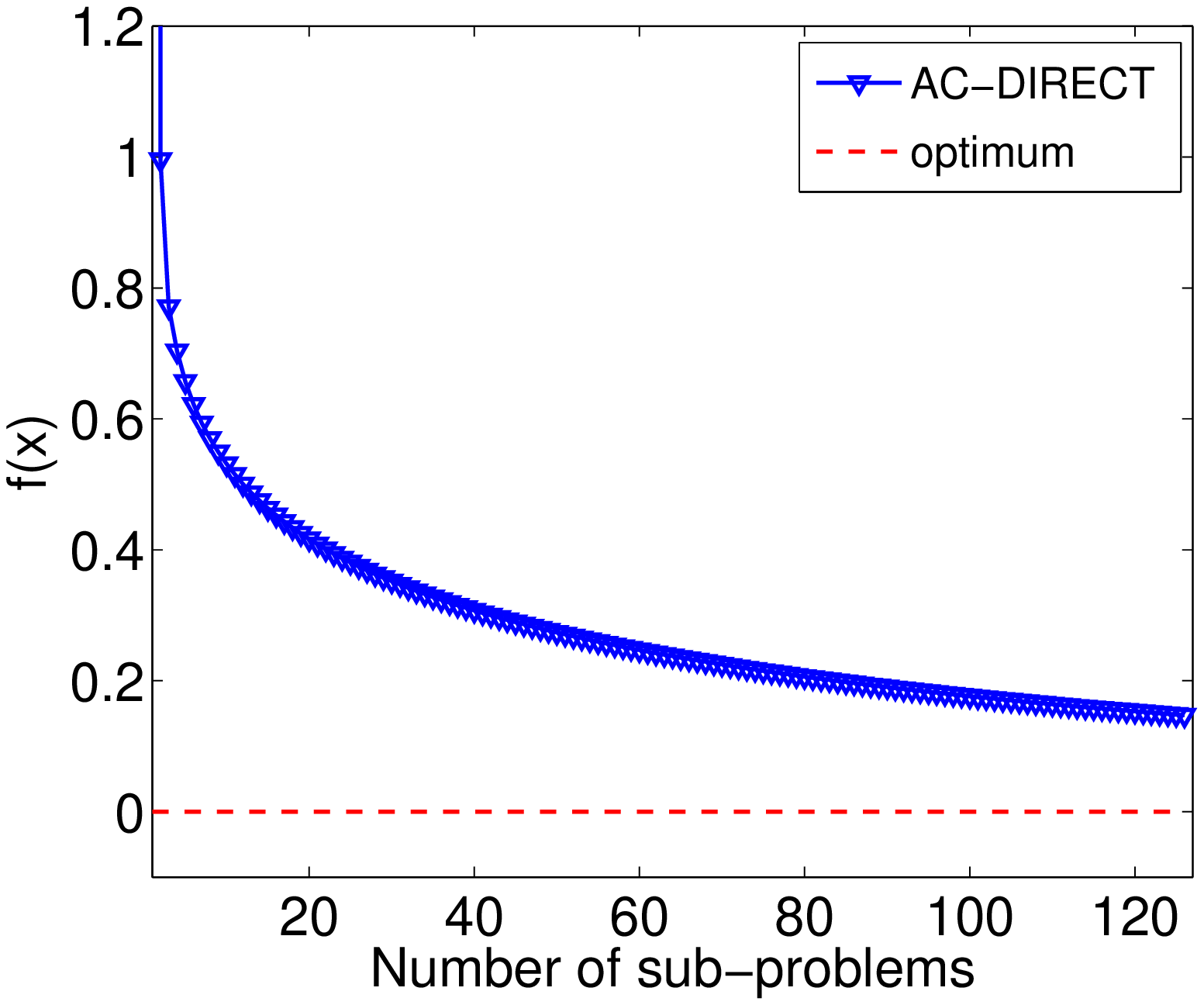}
}
\subfigure[] { \label{fig:b}
\includegraphics[width=0.45\linewidth]{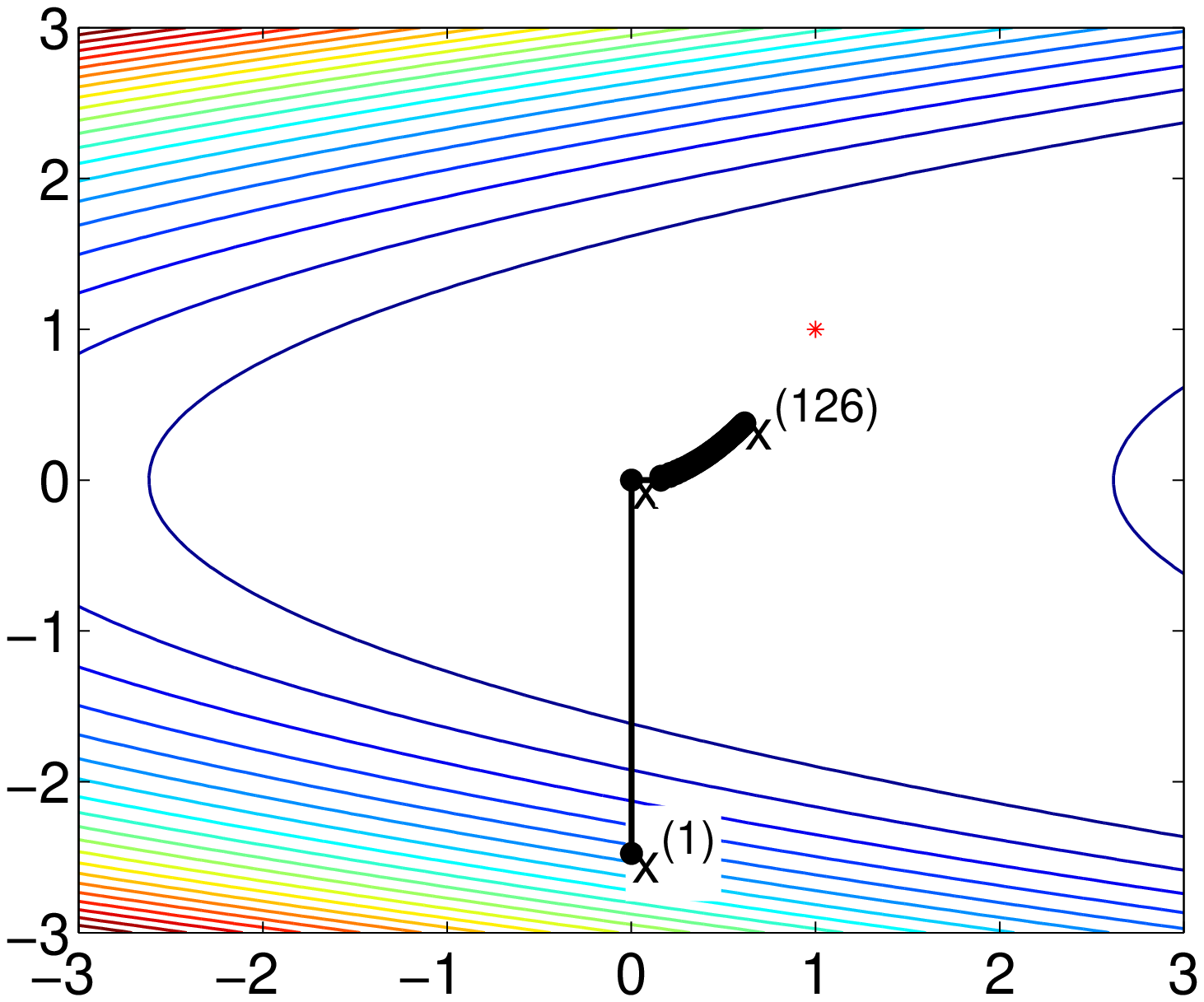}
}
\caption{Simulation results of Rosenbrock function when sequentially applying one-dimensional DIRECT in each sub-problem. Sub-figure (a) illustrates descent of $f(x)$ in every sub-problem. Sub-figure (b) depicts the searching trace on contour map, where the red asterisk is denoted as the global optimal.   }
\end{figure}

Fig.4 shows that one-dimensional DIRECT brings smaller descent to the objective function $f(\mathbf{x})$ as sub-problems go on. The contours around the global optimum are very sparse, and our method only moves forward a trivial step to get closer to the global optimum in each sub-problem. Thus, tremendous numbers of sub-problems are required to approach the global optimum. Although each one-dimensional sub-problem runs efficiently with flat objective function, the huge number of sub-problems also brings inefficiency in solving the original problem. As mentioned in Section 3.2, in such case, when one-dimensional DIRECT appears trivial descent, we switch to block coordinate update. However, the switch remains undesirable results. Therefore, we turn to employ a local optimizer to tackle the problem where the contour is sparse. The employment of local optimizer can potentially help conquer the trivial descent around smooth area. In this paper, SQP is introduced to do the local search. SQP solves a sequence of optimization sub-problems, each of which optimizes a quadratic model of the objective function. If the problem is unconstrained, then the method reduces to Newton's method for finding a point where the gradient of the objective vanishes. SQP can be implemented with fast solving speed in Matlab, thus the incorporation with it as a local optimizer brings little expense to the speed of the proposed algorithm. For the experiment presented in Fig.4, we apply SQP to Rosenbrock function in the $127th$ sub-problem. The results are shown in Fig.5.

\begin{figure}[!hbt]
\centering
\subfigure[] { \label{fig:a}
\includegraphics[width=0.45\linewidth]{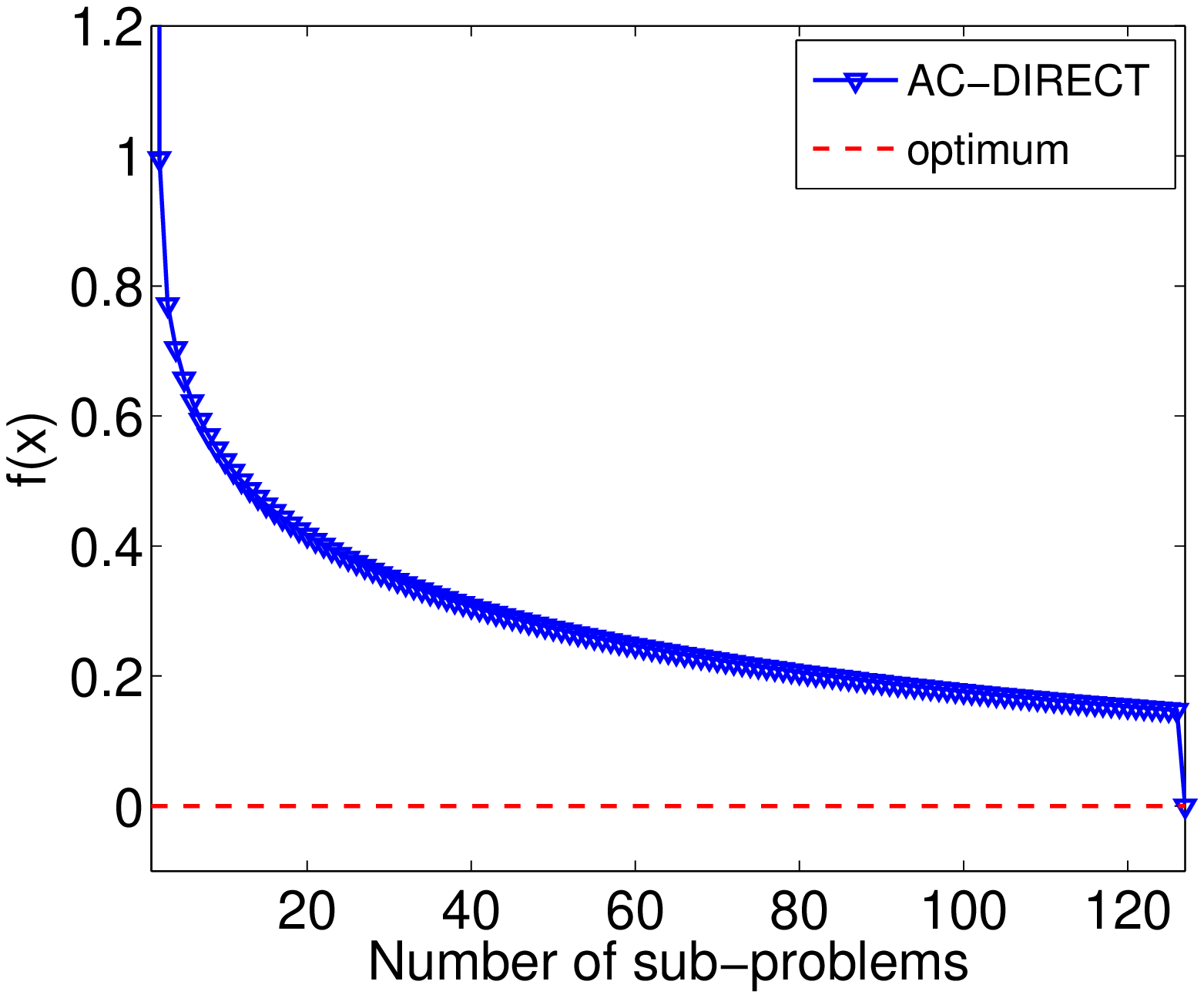}
}
\subfigure[] { \label{fig:b}
\includegraphics[width=0.45\linewidth]{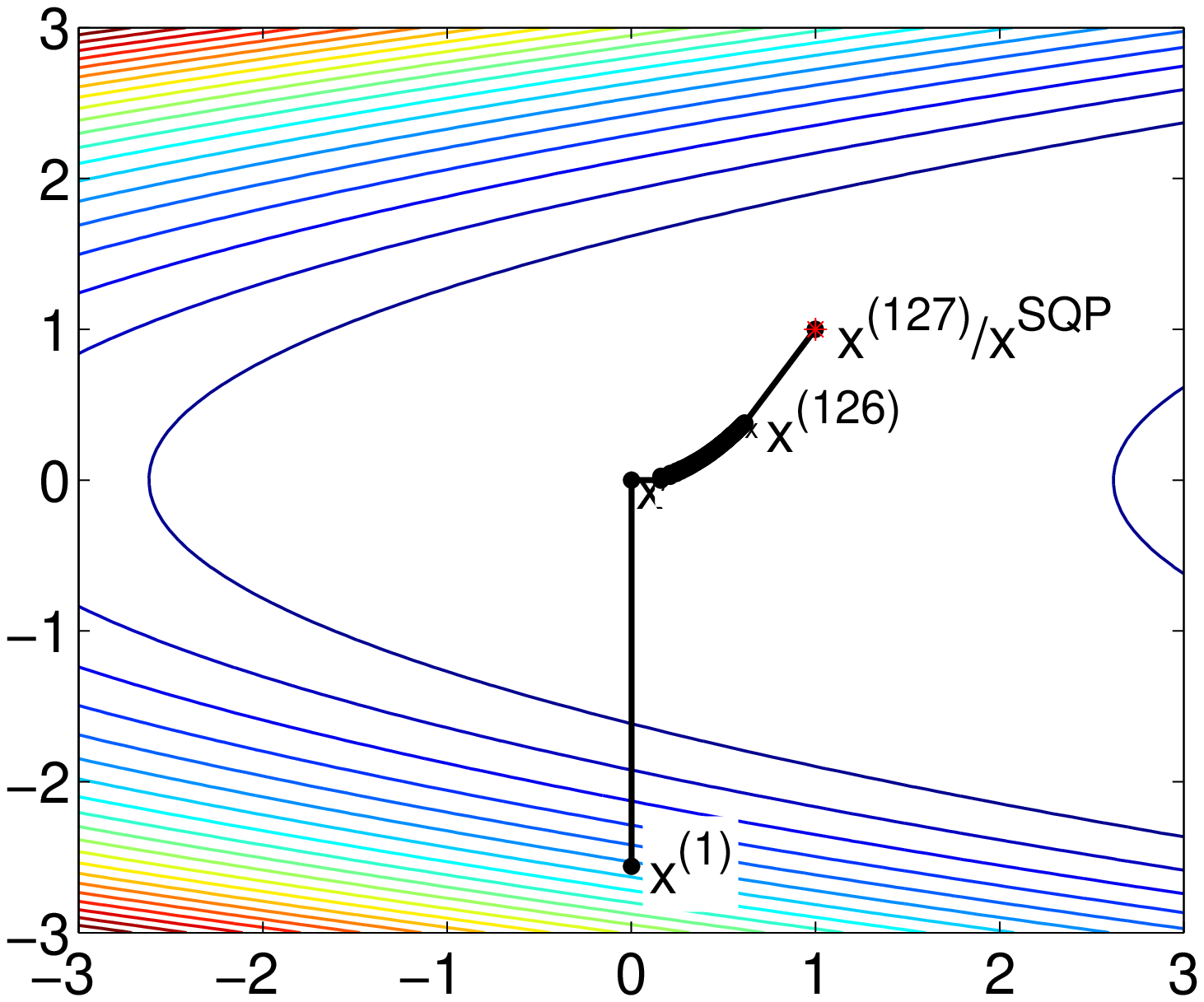}
}
\caption{Based on Fig.4, simulation results of Rosenbrock function with SQP as the local optimizer. }
\end{figure}

Fig.5 demonstrates that the employment of SQP makes the algorithm to quickly approach the global optimum, which presents the efficacy of SQP around smooth area. Integrating the ideas in Section 3.1 and Section 3.2, we establish the ABCD algorithm. The detail is presented in Algorithm 3, where $m_i = 1$ refers to the coordinate uodate and $m_i>1$ represents the block coordinate update.

\begin{algorithm}
\caption{Adaptive Block Coordinate DIRECT algorithm }
\footnotesize{
\KwIn{$f(\mathbf{x}), f^*, \epsilon, \epsilon_1, T_1, m_1, m_2$.}
\KwOut{$f_{\min}$}
Select a feasible starting point $\mathbf{x}_0$, and regard it as the current optimal point$\mathbf{\hat{x}}^*=\mathbf{x}_0$, and the minimal value $\hat{f}^*=f(\mathbf{\hat{x}}^*)$.\\
Choose coordinates $\mathbf{x}(I_1)$ to optimize. \\
The index ${I_1}$ contains $m_1\le n$ elements, which are selected from $\{1,...,n\}$. \\
\While{$|\hat{f}-\hat{f}^*| \leq \epsilon_1$ doesn't happen $T_1$ times in succession and  $|\hat{f}^* - f^*|\geq \epsilon$}
{
Conduct DIRECT only in  the $m_1$ coordinates of $I_1$, and keep the rest fixed to $\mathbf{x}^*(i),i\notin I_1$.\\
 Record the optimal point $\hat{\mathbf{y}}\in \mathbb{R}^{m_1}$ of the current  $m_1$-dimensional DIRECT. \\
Update the current optimal point $\mathbf{x}^*(I_1)=\hat{\mathbf{y}}$ and the minimal value $\hat{f}=f(\mathbf{x}^*)$\\
Reset index $I_1$.
}
\If{$|\hat{f}^* - f^*|\geq \epsilon$}
{
 Conduct local optimizer SQP, and starts from point $\mathbf{\hat{x}}^*$.
}
\While{ $|\hat{f}^* - f^*|\geq \epsilon$}
{
Conduct DIRECT only in the $m_2$ coordinates of $I_2$, and keep the rest fixed to $\mathbf{x}^*(i),i\notin I_2$.\\
Record the optimal point $\hat{\mathbf{y}}\in \mathbb{R}^{m_2}$ of the current $m_2$-dimensional DIRECT. \\
Update the current optimal point $\mathbf{x}^*(I_2)=\hat{\mathbf{y}}$ and the minimal value $\hat{f}=f(\mathbf{x}^*)$\\
Reset index $I_2$.
}
$f_{\min}=\hat{f}^* $ }
\end{algorithm}

\subsection{Algorithm Implementation}
In Algorithm 3, there are many adjustments which can be adaptively done to better tackle different problems. That is  the reason for  why we name it as an adaptive algorithm. In every sub-problem of Algorithm 3, a single coordinate or a block of coordinates can be selected to get optimized. Hence, the size of chosen coordinates and the way of choosing them can be adapted in varied manners. The starting point $\mathbf{x}_0$ is another crucial factor influencing the final results. Moreover, the inserting of local optimizer SQP is also adjustable. These potential modifications can be regarded as the flexibilities of the proposed algorithm. Due to these flexibilities, the proposed ABCD algorithm is capable of  tackling different problems to obtain better solutions.

1) Starting point

The proposed algorithm initially starts from a point $\mathbf{x}_0$, and then only the selected coordinates conduct DIRECT. A good starting point may help accelerate the convergence of $f(\mathbf{x})$ and approach a better solution within less iterations. The starting point $\mathbf{x}_0$ can be generated anywhere in feasible domain. Inspired by SOO, we can introduce the optimistic choosing mechanism to the starting point selection. SOO positively chooses to divide the sub-domain with the lowest objective value. Similarly, we generate a set of points and select the one with the lowest objective value to start the algorithm.

In implementation, we initially divide the domain into $q$ sub-domains and randomly generate points $\hat{\mathbf{x}}_i,i=1,...,q$ in each sub-domain. Next, we choose the sub-domain $j$ with the lowest objective value as the starting point, i.e.,  $j = \{i|\min\{f(\hat{\mathbf{x}}_i), i=1, ..., q\}\}$. Since $\hat{\mathbf{x}}_j$ is the best point among the sampled points, we positively believe that sub-domain $j$ possibly contains comparatively lower objective values in whole. The selection of starting points only requires $q$ function evaluations, whose time expenses can be almost ignored. In fact, the $q$ sub-domains can be selected in different ways, and the value of $q$ can also be adjusted. Furthermore, other methods of selecting the starting point are also acceptable if they take little cost in computing and bring possible improvements to the latter optimization procedures.

2) Size of chosen coordinates

In every sub-problem, the proposed ABCD algorithm focuses on updating a single coordinate or a block of coordinates, while keeping the rest fixed. First, we select a single coordinate to optimize. One-dimensional DIRECT algorithm runs very fast in each sub-problem, but it requires more sub-problems. Meanwhile, for some problems, one-dimensional optimization may ignore the potential of further improvement among blocks of coordinates. Therefore, in Algorithm 3, the size $m$ of chosen coordinates can be adjusted. The size of chosen coordinates, denoted by $m$, is unnecessarily kept the same as optimization goes on, thus we introduce a switch.  This switch is designed to check  further descent of the objective function $f(\mathbf{x})$ with other size of chosen coordinates. We can choose different $m$ and the switch condition $\epsilon_1, T_1$ is also adjustable. In this paper, we firstly set $m_1=1$ to choose only one coordinate to optimize. If the speed appears undesirable, we switch to the block coordinate update with $m_2=2$.

3) The way of coordinates choosing

Once we decide the size of chosen coordinates in each sub-problem. The way of choosing these coordinates also affects the performance. Sequentially optimizing the coordinates is firstly considered as a simple way, where each coordinate has the same chance to get explored and optimized. This mode presets a certain sequence of coordinates to get optimized, and it equals repeating the optimization on the same coordinates after a round of optimizing all the coordinates. For some problems, sequential optimization may brings slow speed, since the objective function $f(\mathbf{x})$ can happen to enjoy slow convergence under the current mode of coordinates choosing. Thus, randomly choosing coordinates to optimize may help get out of the current mode to explore further improvements. The efficacy of randomized BCDM for convex problems are illustrated in \cite{Refrbcdm,Refconvex}. For optimization problems with convex objective function, reference \cite{Refrbcdm} proves that the existence of the minimal iteration for any given mathematical expectation of the objective function. That is to say, for convex problems, desirable results can be expected with enough iterations. Although the minimal iteration can not be specifically deduced for problem (1) where the convexity is unguaranteed, the spirits of randomized BCDM are also applicable in the proposed algorithm.  Based on the test in Fig.3, instead of sequentially selecting two coordinates to optimize, we turn to randomly choose two coordinates to get optimized in each sub-problem.
\begin{figure}[!hbt]
\centering

\includegraphics[width=0.6\linewidth]{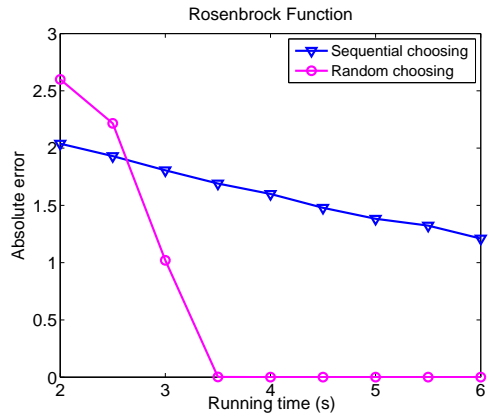}

\caption{Simulation results of Rosenbrock function with ABCD algorithm. `Sequential choosing' and `Random choosing' represent the results of sequentially selecting two coordinates and randomly choosing two coordinates to conduct DIRECT.   }
\end{figure}

Fig.6 illustrates that randomly choosing two coordinates to optimize brings sharper descent compared with sequential choosing. In this paper, we firstly sequentially select one coordinate to optimize, which enables every coordinates share the same chance to update. When we begin the switch, we turn to the mode of random choosing to check further descent of objective function $f(\mathbf{x})$.

4) Local optimizer

In Algorithm 3, when the coordinate update starts to bring trivial descent in each sub-problem, we introduce a local optimizer SQP to accelerate the speed around smooth area. In fact, the switch condition $m_1$ and $T_1$ can also be adjusted. Smaller values of $m_1,T_1$ allows SQP to start the local search at a earlier time in the proposed algorithm, which means less tolerance for cost in CDM. While bigger values of $m_1,T_1$ allow more tolerance in running time spent on CDM. A fixed setting of $m_1,T_1$ may not be able to achieve fairly satisfactory solutions for all test functions, thus the inserting of the local optimizer is not necessary after the CDM and before the BCDM. Sometimes, SQP performs very well when it gets applied to the test function directly. In such case, we can relax the switch condition to introduce SQP earlier or just directly place SQP at the beginning of the algorithm.

\section{Numerical Experiments}
In this section, we conduct numerical experiments to evaluate the proposed algorithm, especially for problems in high dimensions. A series of comparisons among the proposed ABCD algorithm, DIRECT, SOO and MCS are presented. The toolbox of DIRECT, given by Jones, lacks efficiency in the implementation with Matlab. Thus, Bj\"{o}rkman discussed the efficient implementation of DIRECT in Matlab in \cite{Ref19}. To distinguish it from DIRECT, the implementation of DIRECT in \cite{Ref19} is called as glbSolve. In glbSolve, the core idea of identifying the POHs and domain partition still resemble that of DIRECT. For storing the information of each rectangle, instead of using the tree structure in Jones' toolbox, glbSolve uses a straightforward matrix/index-vector technique. It is proved in \cite{Ref19} that glbSolve outperforms Jones' DIRECT toolbox in implementation. Thus, we adopts glbSolve \cite{Ref19} toolbox to implement DIRECT algorithm in numerical experiments, which enables the experimental comparisons more impartial.

Previous researches show that DIRECT can efficiently solve low-dimensional problems, hence we first apply the proposed algorithm to Jones test set, which is regarded as benchmarks for comparing global searching methods in low dimensions. In Section 4.1, we tentatively show the efficacy of coordinate update in DIRECT, and demonstrate that our algorithm maintain superiority in low dimensions. Next, in Section 4.2, we extend Jones test set to Hedar test set to evaluate the performance on more test functions, whose dimensions are adjustable to higher dimensions. Section 4.2 is designed to compare the proposed algorithm with more relevant algorithms,  when tackling different functions in different dimensions. To avoid randomness, the experiments are repeated 5 times for each test function, thus figures in this section are the average values. All the experiments are run on the Windows 7 platform of Matlab R2013a in Core(TM) i7-4790 3.60 GHz, 16GB RAM.

\subsection{Test on the Jones test set }

The Jones test set is regarded as benchmark for comparing global searching methods in low dimensions\cite{Ref1,Ref11,Ref4,Ref8,Ref9,Ref10}. We compare with the original DIRECT to show the efficacy of the coordinate update in the proposed algorithm even for low-dimensional functions.

\begin{table}[!htb]\label{table1}
    \begin{center}
    \caption{Functions of Jones test set.}
     \vskip 0.3cm
    \footnotesize{
    \begin{tabular}{lllll}
        \hline
    Function& Abbreviation  &  Dimension  & Local optimum & Global optimum \\
   \hline
  Shekel 5 & S5 & 4 & 5 & 1\\
  \hline
    Shekel 7 & S7 & 4 & 7 & 1\\
  \hline
    Shekel 10 & S10 & 4 & 10 & 1\\
  \hline
    Hartman 3 & H3 & 3 & 4 & 1\\
  \hline
    Hartman 6 & H6 & 6 & 4 & 1\\
  \hline
    Branin RCOS & BR & 2 & 3 & 3\\
  \hline
    Goldstein and Price & GP & 2 & 4 & 1\\
  \hline
    Six-Hump Camel & C6 & 2 & 6 & 2\\
  \hline
    2D Shubert & SHU & 2 & 760 & 18\\
  \hline
    \end{tabular}}
    \end{center}
\end{table}

In Jones test set, there are 9 problems which are abbreviated as S5, S7, S10, H3, H6,
BR, GP, C6, and SHU. Table 1 shows that the test functions in Jones test set have multiple local optima, and some even have multiple global optima. Similarly, we perform the test as done in \cite{Ref8}. Each algorithm is terminated
when
\begin{equation}
|f_{\min} - f^{*}|\leq \epsilon ,
\end{equation}
where $f_{\min}$ is the minimal value obtained by the tested algorithms. $f^*$ is the global minimum of the test function, and $\epsilon $ is the given accuracy. Since the test functions are low-dimensional ($n \in [2, 6]$), sequentially optimizing the coordinates one by one is capable of obtaining desirable results. Thus, in Algorithm 3, we set $m_1=1$ and local optimizer as well as the switch are ignored. The results are summarized in Table 2. If algorithm fails to stop within 10s CPU running time, we say that it fails in solving the problem within acceptable running time and mark `Fail(t)' in Table 2.
\begin{table}[!htb]\label{table2}
    \begin{center}
    \caption{Running time (s) for reaching accuracy $\epsilon_1$ and $\epsilon_2$, where $\epsilon_1 = 10^{-4}$ and $\epsilon_2=10^{-6}$. `ABCD(1)' means $m_1=1$ in Algorithm 3, where the procedures of the local optimizer SQP and the switch to other size of chosen coordinates are ignored.  }
    \footnotesize{
    \begin{tabular}{llllllllll}
\hline
  & S5 & S7& S10 & H3 & H6 & BR & GP & C6 & SHU \\
\hline
 DIRECT: $\epsilon_1$ & 0.037 & 0.038& 0.039 & 0.037 & 0.055 & 0.029 & 0.031 & 0.031 & Fail(t) \\
\hline
ABCD(1): $\epsilon_1$     & 0.056 & 0.053& 0.046 & 0.037 & 0.069 & 0.029 & 0.049 & 0.032 & 0.034 \\
\hline
 DIRECT: $\epsilon_2$  & 0.066 & 0.061&  0.063 & 0.079 & 0.103 & 0.034 & 0.039 & 0.034 & Fail(t) \\
\hline
ABCD(1): $\epsilon_2$      & 0.051 & 0.057&  0.058 & 0.041 & 0.079 & 0.057 & 0.035 & 0.034 &  0.035 \\
\hline
    \end{tabular} }
    \end{center}
\end{table}

Table 2 illustrates that the proposed algorithm is comparable to DIRECT,  when the accuracy is selected as $\epsilon_1= 10^{-4}$. The proposed ABCD algorithm iteratively selects one coordinate to optimize, and the test function self is low-dimensional. Thus, it makes sense that the proposed algorithm may not be distinctively superior to DIRECT for reaching a relatively low accuracy  in low dimensions. However, when accuracy is restricted to a higher level, say $\epsilon_2 = 10^{-6}$, DIRECT converges slower and slower in the fairly close vicinity of the global optimum. While, the proposed algorithm converges faster that DIRECT for $\epsilon_2=10^{-6}$.

It is worth to mention that, for Shubert function, DIRECT fails both in $\epsilon_1$ and  $\epsilon_2$. Since there are 760 local optima and  18 global optima, DIRECT gets trapped in a local optimum and is unable to get out of it, thus it fails to reach the given accuracy either for $\epsilon_1$ or $\epsilon_2$. While, the proposed algorithm  successfully bypasses the local optima and quickly converges to the global optimum. The proposed ABCD algorithm may cannot guarantee the global optima for each problem within limited budgets. But the results in Table 2 illustrates that it can speed up the convergence, and provide possibilities to bypass the local optima to approach a better solution or even the global optima. Thus, the efficacy of coordinate update in ABCD algorithm is tentatively proven when compared with DIRECT in low dimensions.

\subsection{Test on the Hedar test set}
  To further solidify the experiments, we compare the performance of the proposed ABCD algorithm with DIRECT, SOO and MCS with more test function in this subsection. The Jones test is focusing on low-dimensional problems, on which DIRECT has been reported as a very good algorithm. Our method is comparable and slightly better than DIRECT, showing the adaptiveness of ABCD. In the following, we go to Hedar test, which is largely extended from the Jones test set and contains many high-dimensional test functions. Herdar test set is regarded as benchmarks for global search method \cite{Ref1,Ref28}.

  In low dimensions, DIRECT, SOO and MCS are all capable of obtaining a satisfactory solution within desirable running time. In section 4.1, we have already proved the efficacy of the proposed algorithm for low-dimensional functions. Thus, further experiments are established to present the divergent performance of the tested algorithms for high-dimensional functions. In Hedar test set, there are many test functions which are adjustable to high dimensions. We skip the test functions included in the Jones test set already, and also ignore the test functions in $\mathbb{R}^2$ space. We only consider the test functions with relatively larger scales, say dimension $n> 4$. To investigate the influence of dimension $n$, the selected test functions are required to be changeable in dimensions $n$. Therefore, we obtain 13 test functions from Hedar test set to do the following experiments.

In this subsection, we evaluate the performance via CPU running time $t$ and the absolute error $\epsilon$, which is defined in equation (5). We choose $\epsilon=10^{-4}$ as the target accuracy. Each algorithm terminates when reaching $\epsilon$. Considering the case where algorithm converges in a local optimum or algorithm fails to reach the given accuracy $\epsilon$ within the budget of CPU running time $t$, we terminate the algorithm when any of the following criteria gets satisfied.\\
1. $|\hat{f}-f^*|<\epsilon$,\\
2. $\hat{f}$ fails to decrease $\epsilon_1 = 10^{-6}$ for $\min\{n,6\}$ sub-problems in succession,\\
3. $t > 20s$.

 In many functions, the optimal point is right on the axis origin and the searching domain is axis symmetrical at the same time. However, DIRECT first samples the center point of the initial domain, and then it reaches the global optimum at the first iteration. To guarantee the fairness of the comparisons , in such case, the lower bounds and upper bounds are set asymmetrical. Referring to \cite{Ref1}, the lower bound $\mathbf{L}$ is adjusted to $0.8\mathbf{L}$ and the upper bound $\mathbf{U}$  is changed to $1.2\mathbf{U}$ when the mentioned symmetry happens.

 In this subsection, we exert the proposed algorithm with its flexibilities. We first incorporate the coordinate update, which means sequentially conducting one-dimensional DIRECT on every coordinate. When the switch condition is satisfied, local optimizer SQP is introduced. Then we switch to the block coordinate update, which is presented as $m_2=2$. In the switch, we iteratively select 2 coordinates to optimize in randomness instead of in sequence. The switch condition is set as $T_1 = 3$ and $\epsilon_0 = 10^{-3}$ in this paper. That is to say, when the coordinate update fails to decrease $\epsilon_0=10^{-3}$ for $T_1 = 3$ sub-problems in succession, then we switch to SQP and set $m_2=2$ in the block coordinate update.  When an algorithm performs best in $p$ test functions, we say that the winning ratio of this algorithm is $p/13$. Below, Table 3 summarizes the winning ratios of the tested algorithms for Hedar test set.

\begin{table}[!htb]\label{table3}
    \begin{center}
    \caption{Winning ratios of tested algorithms for reaching accuracy $\epsilon=10^{-4}$. `ABCD(1-2)' represents the results of the proposed Algorithm 3 where $m_1=1$ and $m_2=2$.} \vskip 0.3cm
    \footnotesize{
    \begin{tabular}{ccccc}
\hline
 Algorithm& SOO& MCS & DIRECT& ABCD(1-2) \\
\hline
Winning ratio & 0 & 2/13 & 0 & 11/13   \\
\hline
    \end{tabular}}
    \end{center}
\end{table}

 Table 3 shows that the proposed ABCD algorithm holds the highest wining ratios for Hedar test set. The proposed ABCD algorithm succeeds in reaching the given accuracy $\epsilon=10^{-4}$ within shorter CPU running time for most of the test functions. To further investigate the performance of the tested algorithms. The detailed numerical results for each test function are demonstrated in Table 4 and Table 5, which is a supplyment for Table 4.  $t_{SOO}$, $t_{MCS}$, $t_{DIR}$ and $t_{ABCD(1-2)}$ represent the CPU running time of SOO, MCS, DIRECT and the proposed ABCD algorithm. `Fail(t)' means that algorithm keeps achieving descent as time goes on but fails to reach the given accuracy $\epsilon=10^{-4}$ within the budget of CPU running time, while `Fail(l)' represents that algorithm gets stuck into a local optimum and fails to get out of the current local optimum within the budget of CPU running time. In Table 4 and Table 5, for each test function in each selected  dimension, the best performance among the tested algorithms is printed in boldface and marked with an underline. To distinguish the efficacy of coordinate update and local optimizer SQP, the two columns on the right side of Table 4 and Table 5 are presented for reference. We apply the proposed ABCD algorithm with ignoring SQP to the test function, where the running time is denoted as $t^0_{ABCD(1-2)}$ . We also independently apply SQP to the test function, where the running time is written as $t_{SQP}$.

\begin{table}[!htb]\label{table4}
    \begin{center}
    \caption{Running time (s) for reaching given accuracy $\epsilon=10^{-4}$.}
     \vskip 0.3cm
    \footnotesize{
    \begin{tabular}{llllll|ll}
\hline
 Function & Dim & $t_{SOO}$ & $t_{MCS}$ & $t_{DIR}$  & $t_{ABCD(1-2)}$ & $t^0_{ABCD(1-2)}$ & $t_{SQP}$ \\ \hline
 Ackley& 6  & 1.173 & 0.471 &3.041  &\underline{\textbf{0.097}} &0.097 & Fail(l)\\ \hline
      & 12  & 8.375 & 0.701  &Fail(t) & \underline{\textbf{0.158}}& 0.158& Fail(l)\\ \hline
      & 18  & Fail(t) & Fail(t) &Fail(t) &\underline{\textbf{0.226}} &0.226 & Fail(l)\\ \hline
 Dixon-Price & 6  & Fail(l) & Fail(l) & Fail(l)  & \underline{\textbf{0.545}}& 0.703&Fail(l)\\ \hline
             & 12 & Fail(l) & Fail(l) & Fail(l)  & \underline{\textbf{0.708}} & 1.583 &Fail(l)\\ \hline
             & 18 & Fail(l) & Fail(l)  & Fail(l)  & \underline{\textbf{1.161}}&4.139 & Fail(l)\\ \hline
 Griewank    & 6  & {{0.179}}  & Fail(t)  &Fail(l)   & \underline{\textbf{0.118}}  & 0.118 &Fail(l)\\ \hline
             & 12 & {{0.798 }} & Fail(t)   &Fail(l)   &\underline{\textbf{0.161}} & 0.161 &Fail(l)\\ \hline
             & 18 & {{2.240}}  & Fail(t)   &Fail(l)   & \underline{\textbf{0.213}}& 0.202 &Fail(l)\\ \hline
  Levy       & 6  & 0.212  & 0.433  & 0.075 & \underline{\textbf{0.072}}  &0.071 &Fail(l) \\ \hline
             & 12 & 1.286  & 0.488 &1.351  &\underline{\textbf{0.112}} &0.113 &Fail(l)\\ \hline
             & 18 & 4.185  & 0.498  &Fail(t) &\underline{\textbf{0.159}}&0.159 &Fail(l)\\ \hline
Michalewicz  & 5  & Fail(t) & Fail(t)  &Fail(t) & \underline{\textbf{0.071}} &0.064 &Fail(l)\\ \hline
             & 10 & Fail(t) & Fail(t)  & Fail(t)&  \underline{\textbf{0.077}}&0.076 & Fail(l)\\ \hline
 Powell      & 6  & 0.514 & \underline{{\textbf{0.350}}} & Fail(t)& {0.460} & 1.266&0.367\\ \hline
             & 12 & Fail(t)  & \underline{{\textbf{0.439}}} & Fail(t) & {0.706} &Fail(t) &0.413 \\ \hline
             & 18 & Fail(t)& \underline{{\textbf{0.506}}}  &Fail(t) &  {0.916} &Fail(t) & 0.428 \\ \hline
 Rastrigin   & 6  & 0.641 & 0.467 &0.181  & \underline{\textbf{0.077} } &0.076 & Fail(l)\\ \hline
             & 12 & 4.189 & 0.649 &14.23  &\underline{\textbf{0.118}} & 0.117& Fail(l)\\ \hline
             & 18 & 14.25& 0.841  &Fail(t) & \underline{\textbf{0.161}}&0.162 &Fail(l)\\ \hline
 Rosenbrock  & 6  & {{0.500}} & Fail(l)  &Fail(t)  & \underline{\textbf{0.674}} & 3.188 &Fail(l)\\ \hline
             & 12 & {{3.399}} & Fail(l)  &Fail(t) &\underline{\textbf{0.719}} & 4.995 &Fail(l)\\ \hline
             & 18 &{ {10.84}}& Fail(l)  &Fail(t) & \underline{\textbf{1.202}}  & 8.202 &Fail(l)\\ \hline
Schwefel    & 6  & Fail(t)& Fail(t) &Fail(t) & \underline{\textbf{0.091}} &0.091 &Fail(l)\\ \hline
             & 12 & Fail(t)  & Fail(t) & Fail(t) &\underline{\textbf{0.531}} & 0.167&Fail(l)\\ \hline
             & 18 & Fail(t)& Fail(t)   &Fail(t)& \underline{\textbf{0.593}}& 0.237&Fail(l)\\ \hline
             Sphere      & 6  & 0.246 & 0.250  &0.058 &\underline{\textbf{0.067}} &0.066 &0.321\\ \hline
             & 12 & 1.482 & 0.289  &1.134 &  \underline{\textbf{0.101}} &0.101 &0.321\\ \hline
             & 18 & 4.435 & 0.366  &6.029 & \underline{\textbf{0.129}} & 0.128& 0.321\\ \hline
    \end{tabular}}
    \end{center}
\end{table}

\begin{table}[!htb]\label{table4}
    \begin{center}
    \caption{Supplement for Table 4. Running time (s) for reaching given accuracy $\epsilon=10^{-4}$.}
      \vskip 0.3cm
    \footnotesize{
    \begin{tabular}{llllll|ll}
\hline
 Function & Dim & $t_{SOO}$ & $t_{MCS}$ & $t_{DIR}$  & $t_{ABCD(1-2)}$ & $t^0_{ABCD(1-2)}$ & $t_{SQP}$ \\ \hline
 Sum Square  & 6  & 0.229 & 0.247  &0.038 &  \underline{\textbf{0.059}}& 0.059&0.363\\ \hline
             & 12 & 2.426 & 0.278 &0.109 & \underline{\textbf{0.101}}  & 0.101&0.385\\ \hline
             & 18 & 7.775 & 0.332  &6.029 &\underline{\textbf{0.137}}&0.136 &0.446\\ \hline
  Trid       & 6  & Fail(t) &\underline{\textbf{0.372}}    &1.136  &  1.129 & 1.844&0.348\\ \hline
             & 12 & Fail(t)&\underline{{\textbf{0.455}}}  &Fail(t)   &5.592 &Fail(t) & 0.363 \\ \hline
             & 18 & Fail(t)& \underline{{\textbf{0.521}}} &Fail(t)    &16.86& Fail(t)& 0.366\\ \hline
  Zakharov   & 6  & 0.432 & {{0.348} }&Fail(t)&  \underline{\textbf{0.402}} &0.997 &0.365 \\ \hline
             & 12 & 3.827 & {{0.637} }&Fail(t)  &\underline{\textbf{0.426}} & 5.386&0.357 \\ \hline
             & 18 & {{16.44}}& Fail(t)&Fail(t)  &\underline{\textbf{0.539}}  &15.90 &0.416 \\ \hline

    \end{tabular}}
    \end{center}
\end{table}

 Table 4 and Table 5 illustrate that the proposed ABCD algorithm has predominant advantages in most of the test functions.  DIRECT fails to achieve the given accuracy $\epsilon = 10^{-4}$ in 10 test functions, and SOO and MCS also fail in more than 5 test functions.  Although the proposed algorithm fails to achieve 100\% winning ratio for all the test functions under the current settings $m_1$, $m_2$, $\epsilon_1$ and $T_1$ in Algorithm 3. Referring to no-free-lunch theory, it is normal that an algorithm is unable to perform best under any circumstance. It worths to attention that the proposed algorithm succeeds in every test function for accuracy $\epsilon = 10^{-4}$. Especially, for function Dixon-Price, Michalewicz and Schwefel, MCS, SOO and DIRECT all fail to reach the given accuracy $\epsilon =10^{-4}$, while ABCD algorithm successfully solves it with high speed.

 From the two columns on right side of Table 4 and Table 5, it is obvious that coordinate update is the fundamental mechanism of the proposed ABCD algorithm and the employment of SQP assists the coordinate update to better handle more test functions. It can be seen that, for Ackely, Griewank, Levy, Michalewicz, Rastrigin, Sphere and Sun square, the results are mainly contributed by coordinate update. For function Dixon-Price, Powell, Rosenbrock, Trid and Zakharov, the employment of SQP accelerates the convergence and helps improve the efficiency of the proposed algorithm. If SQP is directly applied to test functions, it fails for most of the test functions. Thus, SQP only functions its advantages with the start points which are obtained by coordinate update. In summary, the introduction of coordinate update and the local optimizer SQP are shown to be useful in combination, since they mutually enables the proposed algorithm better solves more different problems.

 Table 4 and Table 5 show that SOO, MCS and ABCD all successfully reach the give accuracy $\epsilon=10^{-4}$ for function Levy, Rastringin, Sphere and Sum Square in dimensions $n=6, 12, 18$. However, the running time $t$ and the changing tendency of $t$ are different with dimension $n$ increasing. Fig. 7 illustrates the running time $t$ of SOO, MCS and ABCD against dimension $n$.

 \begin{figure}[!hbt]
\centering
\subfigure[] { \label{fig:a}
\includegraphics[width=0.45\linewidth]{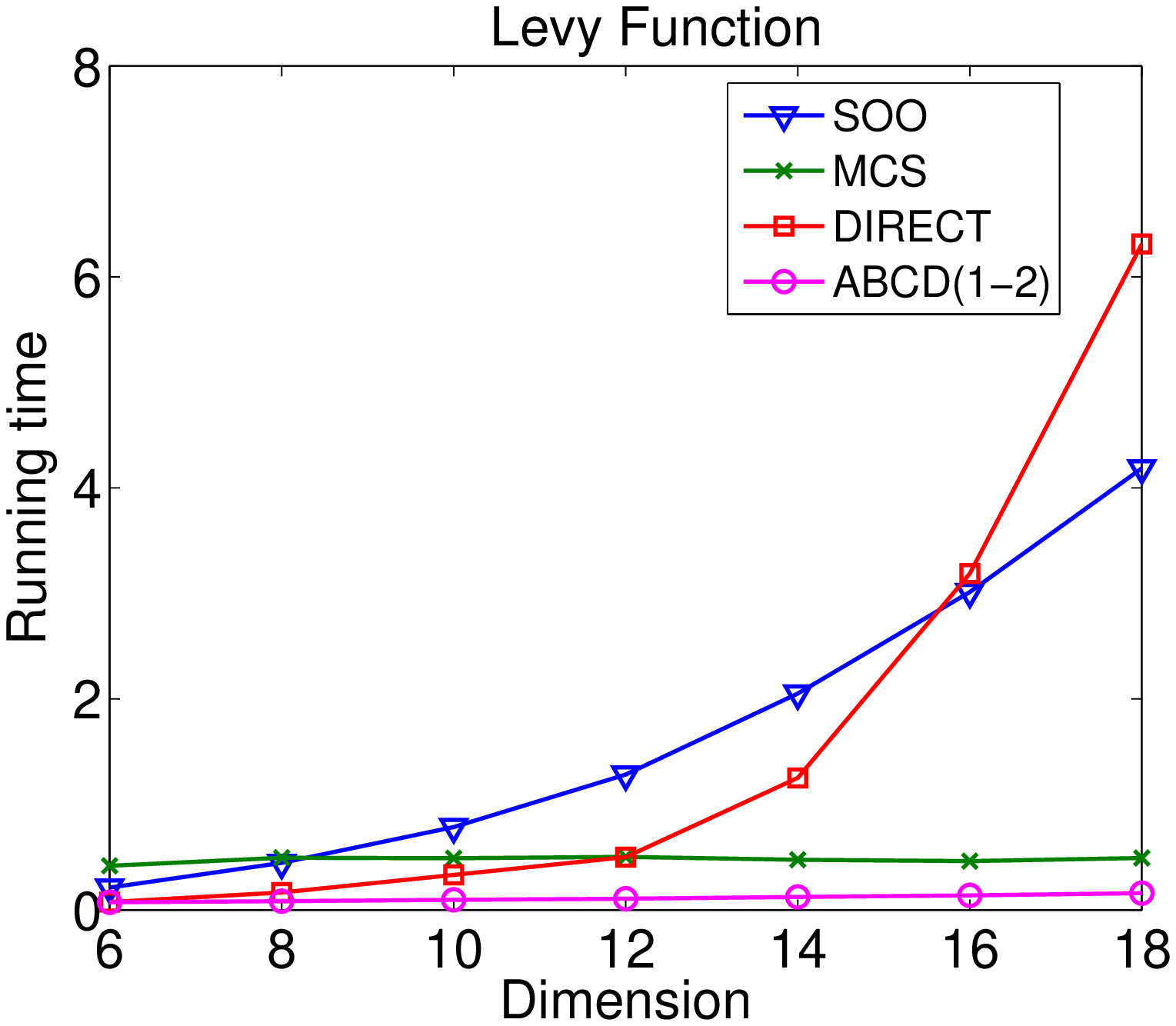}
}
\subfigure[] { \label{fig:b}
\includegraphics[width=0.45\linewidth]{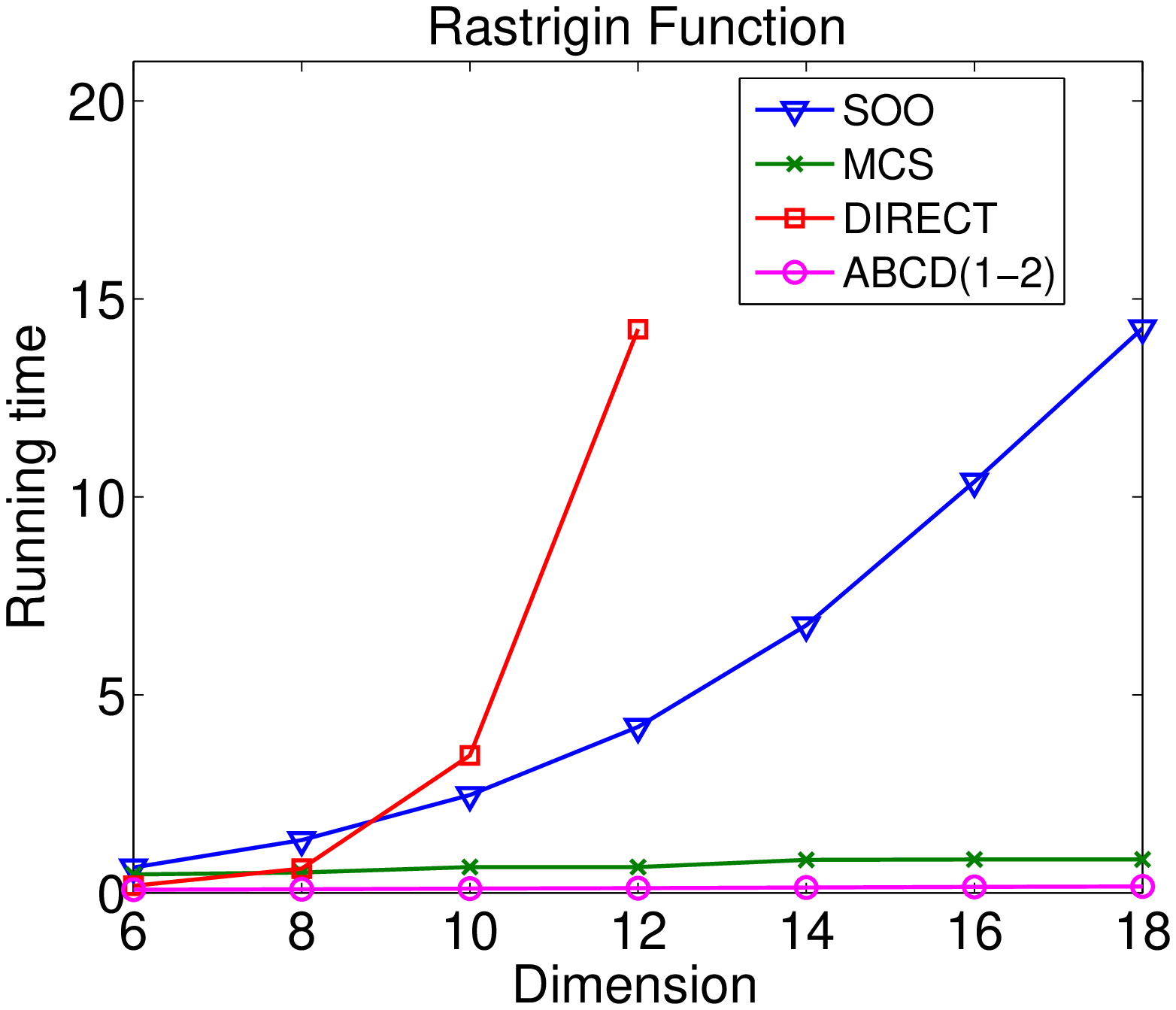}
}
\subfigure[] { \label{fig:c}
\includegraphics[width=0.45\linewidth]{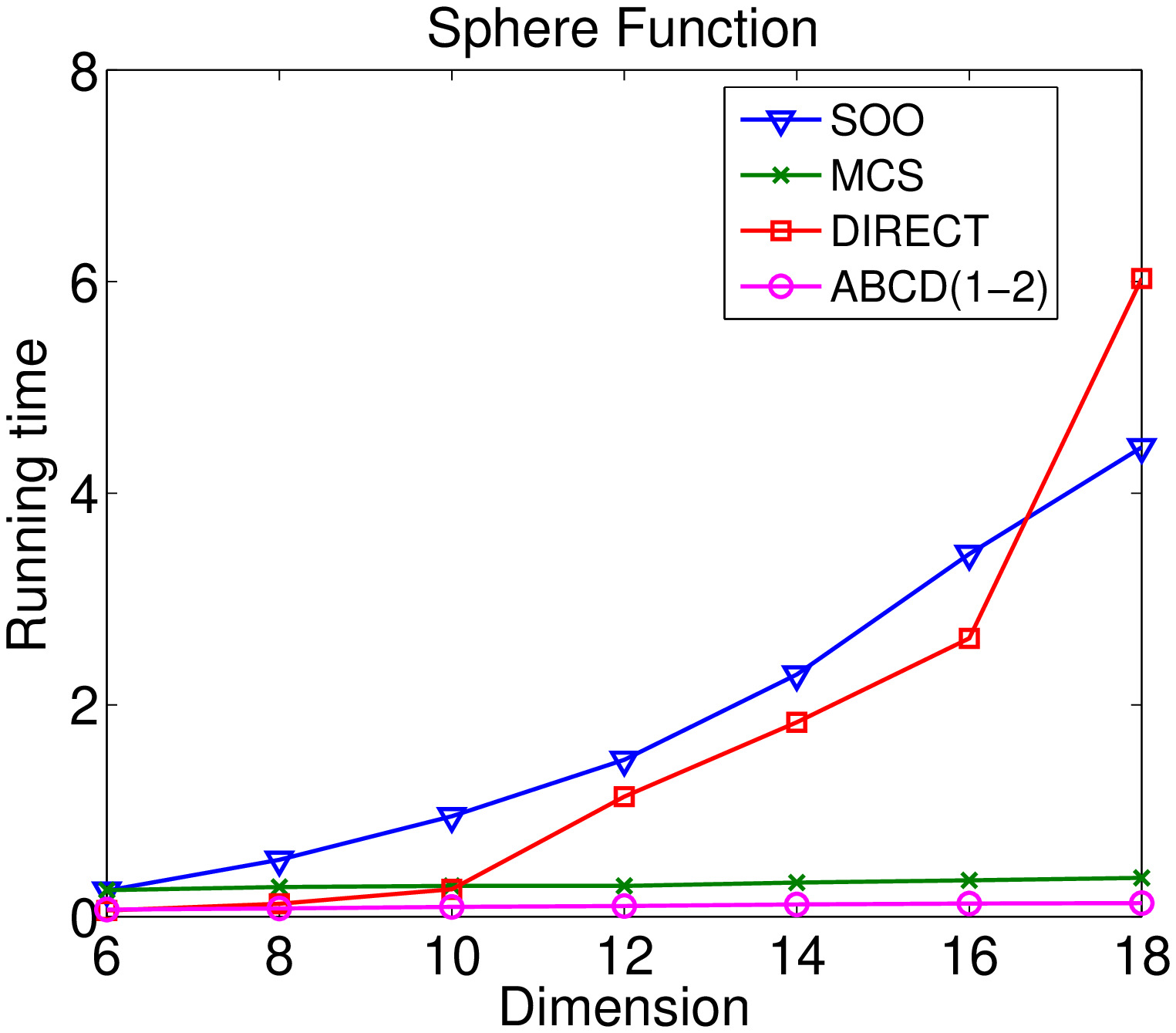}
}
\subfigure[] { \label{fig:d}
\includegraphics[width=0.45\linewidth]{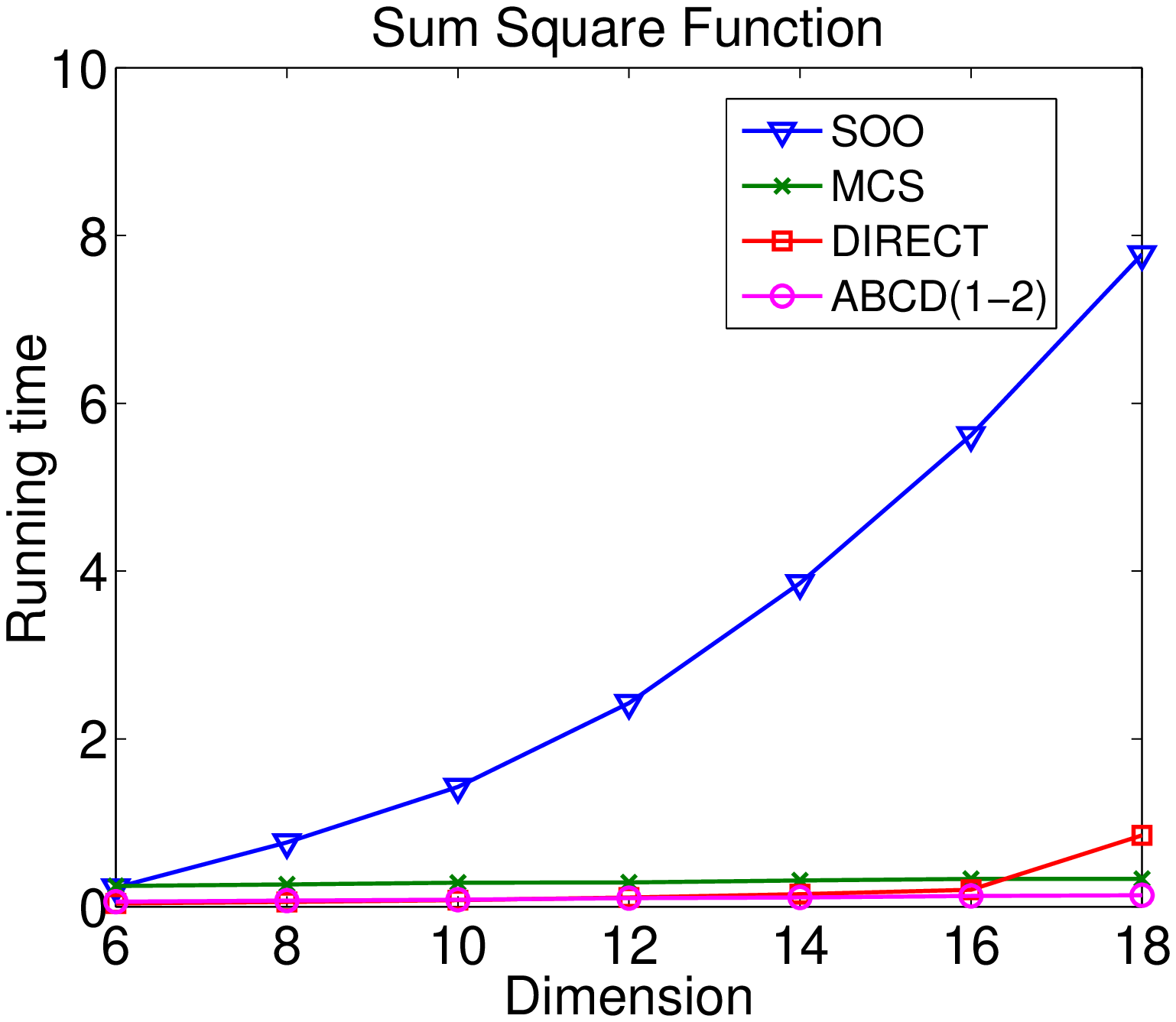}
}
\caption{Sub-test for Levy, Rastrigin, Sphere and Sum Square.  Running time (s) for reaching given accuracy $\epsilon=10^{-4}$. DIRECT fails to achieve accuracy $\epsilon$ within the budget of running time for Rastrigin when dimension $n\geq14$. Thus, the running time of DIRECT for Rastrigin in dimension $n=14, 16, 18$ is omitted.    }
\end{figure}

 Fig. 7 tells that the running time $t_{SOO}$ of SOO increases rapidly with dimension increasing for Levy, Rastrigin, Sphere and Sum Square. Except function Sum square, DIRECT holds the highest slope in Fig.7. Although the running time $t_{MCS}$ of MCS shows a low increasing rate when dimension $n$ increases, the proposed algorithm remains outperforming MCS in speed in every tested dimension. In summary, although SOO, MCS and DIRECT all succeed in solving Levy, Rastrigin, Sphere and Sum Square, the proposed algorithm maintains its advantages in fast speed, and the running time of the proposed algorithm increases slower with dimension increasing.

 For function Powell and  Trid, MCS outperforms our algorithm in the experiment conducted in Table 4 and Table 5. To further investigate the details, we extract Powell function and Trid function individually to present the comparisons, which are illustrated in Table 6. In this sub-test, local optimizer SQP is also evaluated independently to solidify the investigation. In Table 6, $t_{SQP}$ represents the running time $t_{SQP}$  of SQP for reaching accuracy $\epsilon=10^{-4}$, when it is directly applied to the test function.

\begin{table}[!htb]\label{table5}
    \begin{center}
    \caption{Sub-test for Powell and Trid. Running time (s) for reaching given accuracy $\epsilon = 10^{-4}$. }
         \vskip 0.3cm
    \footnotesize{
    \begin{tabular}{ccccc}
\hline
 Function & Dim & $t_{MCS}$ & $t_{SQP}$  & $t_{ABCD(1-2)}$ \\
\hline
 Powell & 6 & \underline{\textbf{0.350}} & 0.367 & 0.460  \\
\hline
        & 12 &0.439  & \underline{\textbf{0.413}} & 0.706  \\
\hline
        & 18 & 0.506 & \underline{\textbf{0.428}} & 0.916 \\
\hline
 Trid & 6 & 0.372 & \underline{\textbf{0.348}} & 1.129 \\
\hline
        & 12 & 0.455 & \underline{\textbf{0.363}} & 5.592 \\
\hline
        & 18 & 0.521 &\underline{\textbf{0.366}} &  16.86 \\
\hline
    \end{tabular}}
    \end{center}
\end{table}

It can be seen from Table 6 that our local optimizer SQP outperforms MCS when it is directly applied to Powell and Trid, let along SQP starts from the point obtained by the coordinate update in the proposed algorithm.  In this experiment, the proposed algorithm first conducts one-dimensional DIRECT to sequentially  optimize the coordinates, then we apply SQP to do the local search when one-dimensional DIRECT is unable to decrease $\epsilon_1=10^{-3}$ in the objective function for $T_1 = 3$ sub-problems in succession. To generalize the experiments on Hedar test set, the switch condition $\epsilon_1$ and $T_1$ are set identical to every test function, even it is sometimes not the optimal choice for a few test functions. It makes sense that the proposed algorithm appears slower than MCS, even though our local optimizer SQP actually outperforms MCS when it is directly applied to Powell and Trid. If we relax the switch condition to conduct SQP earlier, or we replace local optimizer SQP at the beginning, the proposed algorithm still outperforms MCS for Powell and Trid. If we set SQP to begin the proposed algorithm, SQP directly reaches the given accuracy $\epsilon$, the following procedures of  the proposed algorithm are saved. In such case, SQP can be regarded as a degeneration of the proposed algorithm to some extent. This characteristic supports the fact that our algorithm is adaptive and it hold more flexibilities.

\section{Conclusion}
Adaptive block coordinate DIRECT algorithm, presented in this paper, incorporates the strategy of coordinate update to achieve high speed in high dimensions, while DIRECT is shown to loose efficiency and accuracy. Instead of constantly repeating the procedures of sampling and dividing on the whole domain, we iteratively select only one coordinate to do the optimization and keep the rest fixed. Coordinate update maintains the efficiency of DIRECT in low dimensions, but it lacks checking further improvements with other size of chosen coordinates. Thus, we adaptively switch to block coordinate update to explore further improvement with larger size of chosen coordinates. Besides, coordinate update may possible bring trivial descent in each sub-problem when the objective function is very flat, hence local optimizer SQP is employed to accelerate the convergence around smooth area. Furthermore, in the proposed algorithm, the starting point, the size of chosen coordinates, the way of selecting coordinates as well as the inserting of SQP can be adaptively adjusted. These alternative adaptations make the proposed ABCD algorithm more flexible.

\end{document}